\documentclass[11pt]{article}
\usepackage[letterpaper,margin=1in]{geometry}
\usepackage[bf,sf]{titlesec}
\usepackage{algorithm}
\usepackage{algpseudocode}
\usepackage{amsmath}
\usepackage{amssymb}
\usepackage{comment}
\usepackage{fancyvrb}
\usepackage{fixltx2e}
\usepackage{graphicx}
\usepackage{hyperref}
\usepackage{multirow}
\usepackage{pdflscape}
\usepackage{setspace}
\usepackage{tikz}
\usepackage{titling}

\MakeRobust{\Call}

\definecolor{lightgray}{gray}{0.9}

\newcommand*\Let[2]{\State #1 $\gets$ #2}

\makeatletter
\def\@xfootnote[#1]{%
  \protected@xdef\@thefnmark{#1}%
  \@footnotemark\@footnotetext}
\makeatother

\DeclareMathOperator*{\argmin}{argmin}
\DeclareMathOperator*{\minimize}{minimize}
\DeclareMathOperator*{\subjectto}{subject\;to}
\DeclareMathOperator*{\vect}{vec}

\DeclareMathOperator{\prox}{prox}

\title{\textsf{\textbf{Convex programming with fast proximal and linear operators}}}
\author{Matt Wytock, Po-Wei Wang, J. Zico Kolter}

\begin{document}

\maketitle

\begin{abstract}

We present Epsilon, a system for general convex programming using fast linear
and proximal operators. As with existing convex programming frameworks,
users specify convex optimization problems using a natural grammar
for mathematical expressions, composing functions in a way that is guaranteed to
be convex by the rules of disciplined convex programming. Given such an input,
the Epsilon compiler transforms the optimization problem into a mathematically
equivalent form consisting only of functions with efficient proximal operators---an
intermediate representation we refer to as \emph{prox-affine form}. By
reducing problems to this form, Epsilon enables solving general convex problems
using a large library of fast proximal and linear operators; numerical examples
on many popular problems from statistics and machine learning show that this often
improves  running times by an order of magnitude or more vs. existing approaches
based on conic solvers.
\end{abstract}

\section{Introduction}

In the field of convex optimization there has existed a fundamental dichotomy
between ``general purpose'' and ``specialized'' solvers.  Many robust solvers
exist that can solve very general, broad classes of optimization problems
(e.g. semidefinite programs, which can express a very wide range of optimization
tasks), and there are modeling frameworks (e.g. CVX \cite{grant2008cvx}, YALMIP
\cite{lofberg2004yalmip})
that can quickly express complex optimization problems by
reduction to one of these forms.  However, these methods typically only scale to
relatively small-sized problems and largely due to these scaling
problems, there has also been a great deal of work on special-purpose solvers
for very particular optimization problems of interest (e.g. a particular loss
function plus regularization, a particular form of sparse constraints, etc).
But these approaches are very much intended for a particular form of problem, and
usually require significant human effort to adapt to even slightly different
formulations.

In this paper, we present Epsilon, a software architecture and system that
makes a significant step toward bridging this gap.  Like the above-mentioned
modeling frameworks, Epsilon takes as input a convex problem specified in terms
of disciplined convex programming (DCP), a specification of convex problems that can
easily and intuitively express complex objectives and constraints.  However,
unlike current DCP systems which transform the problem into a standard conic
form and then pass to a cone solver, Epsilon transforms the problem into a form
we call \emph{prox-affine}: a sum of ``prox-friendly'' functions (i.e.,
functions that have efficient proximal operators) composed with affine
transformations. These functions include the cone projections sufficient to
solve existing cone problems, but prox-affine form is a much richer
representation also including a wide range of other convex functions
with efficient proximal operators.

The main advantage of the prox-affine form is that it maintains significant
structure from the original problem that is lost in traditional cone
transformations---the
Epsilon solver exploits this structure by directly implementing a large library
of proximal operators and applying operator splitting techniques. In particular,
we develop an approach for solving problems in prox-affine form based on the
alternating direction of multipliers (ADMM) \cite{boyd2011distributed}, but
several alternative approaches are possible as well. Typically, the
computational time for this
algorithm (as with most iterative numerical algorithms) is dominated by the
evaluation of linear operators and here we integrate a large library of
\emph{linear} operators which can often be far more efficient than sparse or
dense matrices. In total, the resulting Epsilon system can solve a wide range of
optimization problems an order of magnitude (or more) that existing approaches
to general convex programming---we provide several examples of popular
problems from statistics and machine learning in Section \ref{sec-examples}.

The main contributions of this paper are:
\begin{enumerate}
\item We develop the prox-affine form for general convex programs and
provide an algorithm for converting general disciplined convex programs into
prox-affine form. Using the terminology of
programming language compilers, the prox-affine form can be thought of as an
intermediate representation (IR) for the convex problem, and multiple
transformations, which make the resulting form easier to solve, can be applied
analogous to the multiple passes of an optimizing compiler.

\item We present an efficient algorithm for solving the resulting transformed
problem, based upon ADMM, which solves the problem using a large library of
linear and proximal operators.

\item We present efficient proximal algorithms for handling a number of
non-linear operators arising from the prox-affine form.  Several of these
operators were also presented in the literature, but we include some additional
cases as well.

\item We develop an efficient system for handling the resulting linear
  operators, which allows us to express complex operators such as
Kronecker products (that frequently arise in these settings) without resorting to
explicit instantiation of sparse/dense matrices.

\item We release an open source version of the Epsilon system, including
  benchmarks for many common optimization problems arising in the machine
  learning and statistics literature. In all cases, we show that Epsilon is
  substantially faster than existing approaches.
\end{enumerate}

\section{Background}
\label{sec-background}

Our work builds heavily upon two main lines of work: 1) disciplined convex
programming (DCP) \cite{grant2006disciplined} and systems built upon
these methods; 2) operating splitting techniques, specifically the alternating
direction method of multipliers (ADMM) \cite{boyd2011distributed}. In this
section we briefly review historical and recent work in these areas, highlighting
how it relates to Epsilon.

\subsection{Disciplined convex programming}
\label{sec-background-dcp}

Disciplined convex programming frameworks such as CVX \cite{grant2008cvx}, CVXPY
\cite{diamond2015cvxpy} and Convex.jl \cite{udell2014convex} allow convex programs to be
specified using a natural programming syntax.  For example, the lasso problem
($\ell_1$-regularized least squares)
\begin{equation}
\minimize_x \|A x - b\|_2^2 + \lambda \|x\|_1
\end{equation}
with optimization variables $x$, problem data $A$, $b$, and regularization
parameter $\lambda$, can be expressed simply as (here using CVXPY,
which our system uses directly):

\begin{verbatim}
x = Variable(n)
prob = Problem(Minimize(sum_squares(A*x - b) + lam*norm1(x))
\end{verbatim}
These libraries and other similar modeling frameworks (e.g. the YALMIP library
\cite{lofberg2004yalmip}) have substantially lowered the barrier to quickly prototyping complex
convex programs without the need to manually convert them to a form that can be
fed directly into a numerical solver (e.g. a standard cone form).

The innovation of DCP methods in particular was a separation between the \emph
{verification} that a program is convex according to DCP ruleset, and the actual
transformation to conic form (and also, naturally, a separation between these
elements and the actual numerical solution).  In particular, DCP methods
represent expressions in the optimization problem (the objective term and all
constraint terms) as abstract syntax trees (ASTs).  Then, a set of composition
rules can be applied to determine whether the problem is in fact a convex one.
For example, suppose $f = h(g_1(x),\ldots,g(k(x))$, for $h : \mathbb{R}^k \to
\mathbb{R}$, $g_i : \mathbb{R}^n \to \mathbb{R}$, and  one of the following
holds for each $i = 1,\ldots,k$:
\begin{itemize}
\item $g_i$ is convex and $h$ is nondecreasing in argument $i$
\item $g_i$ is concave and $h$ is nonincreasing in argument $i$
\item $g_i$ is affine,
\end{itemize}
then $f$ is convex (see e.g. \cite[\S 3.2.4]{boyd2004convex} for
details). Importantly the DCP ruleset is sufficient but \emph{not necessary} for
a problem to be convex: for instance, the ``log-sum-exp'' function is convex,
but is a composition of a concave monotonic and convex function, which does not
imply convexity using these rules.  However, log-sum-exp \emph{can} be
represented as a separate atomic function that is convex (and monotonic in its
arguments).  In practice, most convex problems can be written using the DCP
ruleset with a relatively small set of specially-defined functions.

After convexity of the problem has been verified by the DCP rules, traditionally
the next step in the DCP framework is to convert problems into a standard conic
form.  This is accomplished through epigraph transformations: in addition to
expressing the DCP properties of each function in the DCP collection, the
implementation of DCP function have a representation of the function as the
solution to a linear cone program.  For example, the $\ell_1$-norm can
be expressed as
\begin{equation}
\|x\|_1 \equiv \minimize_t \; 1^T t, \;\; \subjectto \; -t \leq x \leq t.
\end{equation}
Thus, when the transformation step encounters an $\ell_1$-norm, it can be
replaced with the linear cone problem above by simply introducing the $t$
variable, modifying the expression to be that of the objective and adding the
constraints above.

By applying these transformation repeatedly, the entire
problem is reduced to a single linear cone problem which can then be solved by
standard conic solvers. Specifically, the standard problem form used by
disciplined convex programming frameworks is
\begin{equation}
  \begin{split}
    \minimize_x \;\; & c^Tx \\
    \subjectto \;\; & Ax = b \\
    \;\; & x \in \mathcal{K}
  \end{split}
\end{equation}
where $\mathcal{K}$ is the cross product of several cones (e.g. the nonnegative
orthant, second-order and semidefinite cones).  DCP frameworks directly output
problems in this form, typically specifying the $A$ matrix in sparse form.
These problems can then be solved by a wide variety of different approaches, one
of the most common (at least for cases where the cones include the
second order or semidefinite cones) being primal-dual interior point methods
\cite{potra2000interior}. Recent work develops first order methods that enable
scaling to larger problem sizes---the splitting conic solver (SCS)
\cite{o2013operator} is the main comparison for Epsilon in Section
\ref{sec-examples}.


Fundamentally, Epsilon, differs from traditional DCP approaches only in the
intermediate representation that is passed to the solver. We use the exact same
DCP library and ruleset (indeed, our implementation takes problems specified
directly as CVXPY problems), but instead of transforming these problems to
conic form, we apply a different set of transformations which convert the
problem to a higher-level form (the prox-affine form we will discuss shortly);
this intermediate representation maintains a great deal more of the original
problem structure which allows for more efficient direct solutions.

\subsection{Operator splitting}
\label{sec-background-prox}

Operator splitting techniques, such as the alternating direction method of
multipliers (ADMM), have seen a surge of interest in recent years. Essentially,
these methods provide an approach to solving problems with composite objective
\begin{equation}
  \minimize_x \;\; f(x) + g(x)
\end{equation}
by iteratively minimizing each function separately.  A general review of
operator splitting algorithms is given in \cite{ryu2015primer}, and ADMM and
proximal algorithms are reviewed in
\cite{boyd2011distributed,parikh2013proximal} respectively.

Essentially, the main computation of an operator splitting approach reduces to
several application of \emph{proximal operators}: given a function $f :
\mathbb{R}^n \to \mathbb{R}$, the proximal operator is defined as
 \begin{equation}
  \prox_f(v) = \argmin_x \left( f(x) + (1/2)\|x - v\|_2^2 \right),
\end{equation}
where and $\| \cdot \|_2$ denotes the usual Euclidean norm. Conceptually, the
proximal operator can be seen as generalizing set projection to
functions---given an input $v \in \mathbb{R}^n$ we find a point $x$ that is
close to $v$  but also makes $f$ small. We recover the projection operator for a
set $\mathcal{C}$ with the indicator function
\begin{equation}
  I_{\mathcal{C}}(x) = \left\{ \begin{array}{ll}
      0 & x \in \mathcal{C} \\
      \infty & x \notin \mathcal{C},
  \end{array} \right.
\end{equation}
since
\begin{equation}
  \prox_{I_\mathcal{C}}(v) = \argmin_{x \in C} \|x - v\|_2 = \Pi_\mathcal{C}(v).
\end{equation}

Importantly, the key aspect of proximal operators is that they can often be
solved in closed form (or something virtually as efficient) for a wide variety
of functions.  We will detail a large library of these ``prox-friendly''
functions (e.g. functions that have an efficient proximal operator) in Section
\ref{sec-proximal-operators}, and describe the general methods for solving them
efficiently. Epsilon relies crucially on such operators, as it reduces general
optimization problems precisely to a sum of these prox-friendly functions.

\section{The Epsilon system}

\begin{figure}
  \centering
  \includegraphics[scale=0.4]{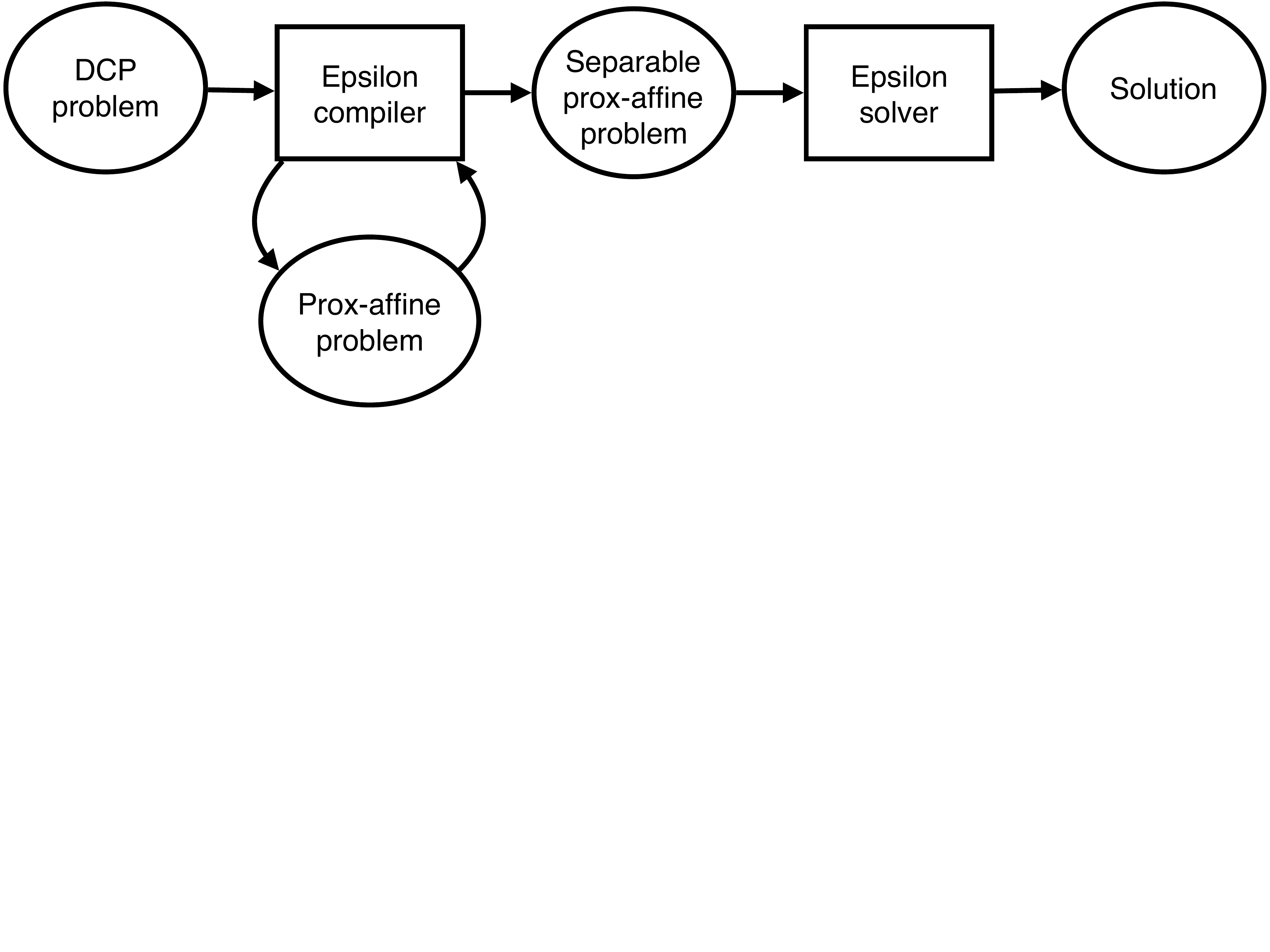}
  \vspace{-2.4in}
  \caption{The Epsilon system: the compiler transforms a
    DCP-valid input problem into separable prox-affine problem to be solved by
    the solver. The prox-affine problem is also used as an intermediate
    representation internal to the compiler.}
  \label{fig-epsilon-system}
\end{figure}

The main idea behind Epsilon is that instead of reducing general convex problems to conic
form, we reduce them directly to a sum of prox-friendly functions composed with
affine transformations and employ a large library of fast proximal and linear
operator implementations to solve problems directly in this form.  We call this
intermediate representation the \emph{prox-affine} form, and problems in this
form can be directly mapped to a set of efficient proximal and linear operators
available to the solver.

There are two main components of Epsilon: 1) the \emph{compiler}, which
transforms a DCP representation into a prox-affine form (and eventually a
separable prox-affine form, to be discussed shortly), by a series of passes over
the AST corresponding to the original problem; and 2) the \emph{solver}, which
solves the resulting problem using the fast implementation of these proximal and
linear operators.  An overview of the system is shown in Figure \ref
{fig-epsilon-system}.  This section describes each of these two components in
detail.

\subsection{The prox-affine form}


The internal representation used by the compiler, as well as the input to the
solver, is a convex optimization problem in prox-affine form
\begin{equation}
  \minimize_x \;\; \sum_{i=1}^N f_i(H_i(x)), \\
\end{equation}
consisting of the sum of ``prox-friendly'' functions $f_1,\ldots,f_N$
composed with a set of affine transformations $H_1, \ldots, H_n$. In Epsilon,
each $f_i \circ H_i$ is \emph{atomic} meaning that it has a concrete numerical
routine in the solver which implements the proximal operator for this function
directly. The affine transformations are implemented with a library of linear
operators which includes matrices (either sparse or dense), as well as more
complex linear transformations like Kronecker products not easily
represented as a single matrix, special cases like diagonal or scalar matrices
and complex chains of multiple such linear transformations. Crucially, not
every proximal operator can be combined with every linear operator, but it is
the job of the Epsilon compiler, described shortly, to ensure that the problem
is transformed to one where the compositions of proximal and linear operators
have an available implementation.  For example, very few proximal operators
support composition with a general dense matrix (the sum-of-squares and subspace
equality constraint being some of the only instances), but several can support
composition with a diagonal matrix. Representing these distinctions in the
Epsilon compiler is critical to deriving efficient proximal updates for the
final optimization problems.

As an example of prox-affine form, the linear cone problem which forms
the basis for existing disciplined convex programming systems (see Section
\ref{sec-background-dcp}),
\begin{equation}
  \begin{split}
    \minimize_x \;\; & c^Tx \\
    \subjectto \;\; & Ax = b \\
    \;\; & x \in \mathcal{K},
  \end{split}
\end{equation}
can be represented in prox-affine form as
\begin{equation}
  \minimize_x \;\; c^Tx + I_{0}(Ax - b) + I_{\mathcal{K}}(x),
\end{equation}
with $f_1$ being the identity function, $H_1$ being
inner product with $c$, $f_2$ the indicator of the zero set, $H_2$ being the
affine transformation $Ax - b$, $f_3$ the indicator of the cone $\mathcal{K}$ and
$H_2$ the identity.  Each of these functions is atomic and therefore has an
efficient proximal operator; for example the proximal operator of $I_0(Ax - b)$
is simply the projection onto this subspace, and the proximal operator for
$I_{\mathcal{K}}$ is the cone projection, and the proximal operator for $c^Tx$ is
simply $v - c$ (in fact, this term can be merged with one or both of the other
terms and thus only two proximal operators are necessary). As the linear cone
problem is thus a special case of prox-affine form, Epsilon enjoys the same
generality as existing DCP systems.



In order to apply the operator splitting algorithm, we also define the
\emph{separable} prox-affine form
\begin{equation}
  \begin{split}
    \minimize_{x_1, \ldots, x_N} \;\; & \sum_{i=1}^N f_i(H_i(x_i)) \\
    \subjectto \;\; & \sum_{i=1}^N A_i(x_i) = b.
  \end{split}
\end{equation}
which has separable objective and explicit additional linear equality
constraints (which are required in order to guarantee that the objective can be
made separable while remaining equivalent to the original problem). As above,
the affine transformations $A_i$ are implemented by the linear operator library
and can thus be represented with matrices or Kronecker products, as well as
simpler forms such diagonal or scalar matrices.  The latter are especially
common in the separable form because they are often introduced for
representing the consensus constraint that two variables be equal (e.g. $A_i =
I$ or $A_i = -I$). Mathematically, there is little difference between the
separable and non-separable forms; but computationally the separable
form allows for direct application of the ADMM-based operator splitting
algorithm. The separable prox-affine form can thus be mapped directly to a
sequence of proximal and linear operator evaluations employed by the operator
splitting algorithm.



The solution methods we present shortly for problems in separable prox-affine
form will ultimately reduce to solving generalized proximal operators of the
form
\begin{equation}
\prox_{f \circ H,A}(v) = \argmin_{x} \lambda f(H(x)) + (1/2)\|A(x) - v\|_2^2.
\end{equation}
For most functions $f$ and affine transformations $H$ and $A$, this function will
not have a simple closed-form solution, even if a simple proximal operator
exists for $f$.  The three important exceptions to this are when: 1) $f$ is the
null function, 2) $f$ is the indicator of the zero cone, or 3) $f$ is a
sum-of-squares function; in all these cases, the solution boils down to a linear
least-squares problem.  However, for certain linear operators (namely, scalar
or diagonal transformations), there \emph{are} many cases where the generalized
proximal operator has a straightforward solution.  A crucial element of the
Epsilon compiler transformations is to produce a prox-affine form where the
combination of $f$, $H$, and $A$ results in a solvable proximal operator; below,
when we list the proximal operators supported by Epsilon, we will also specify
which compositions with linear operators are valid.

\subsection{Conversion to prox-affine form}
\label{sec-compiler-conversion}

\begin{figure}
  \centering
  \includegraphics[scale=0.8]{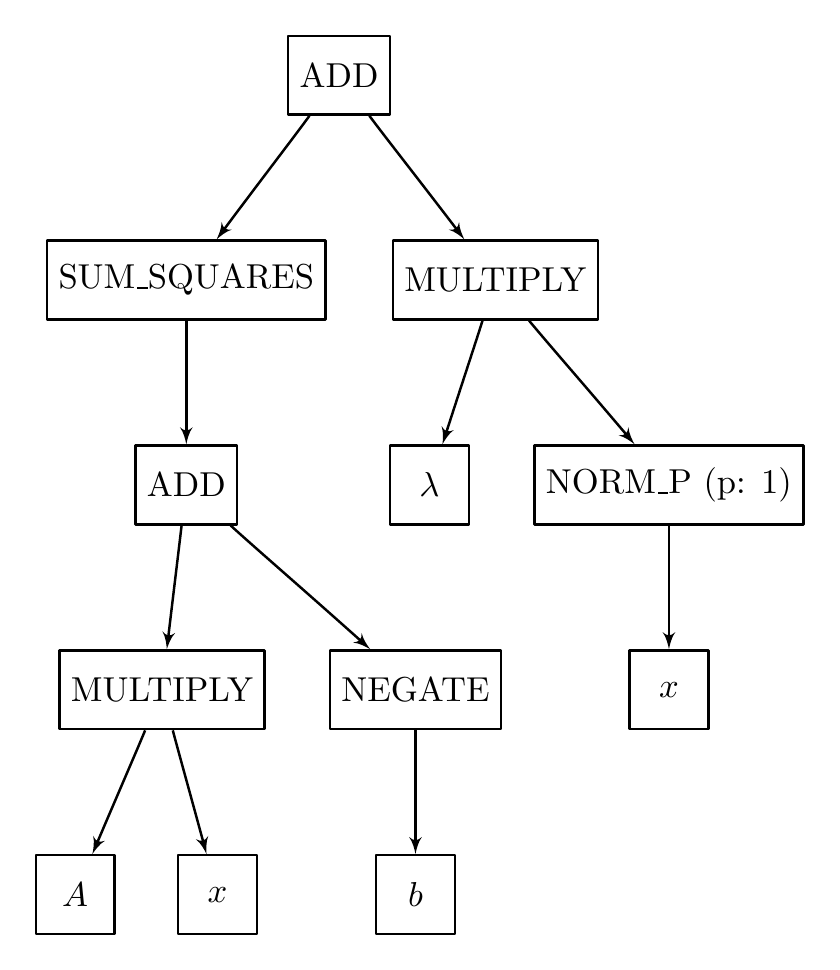}
  \includegraphics[scale=0.8]{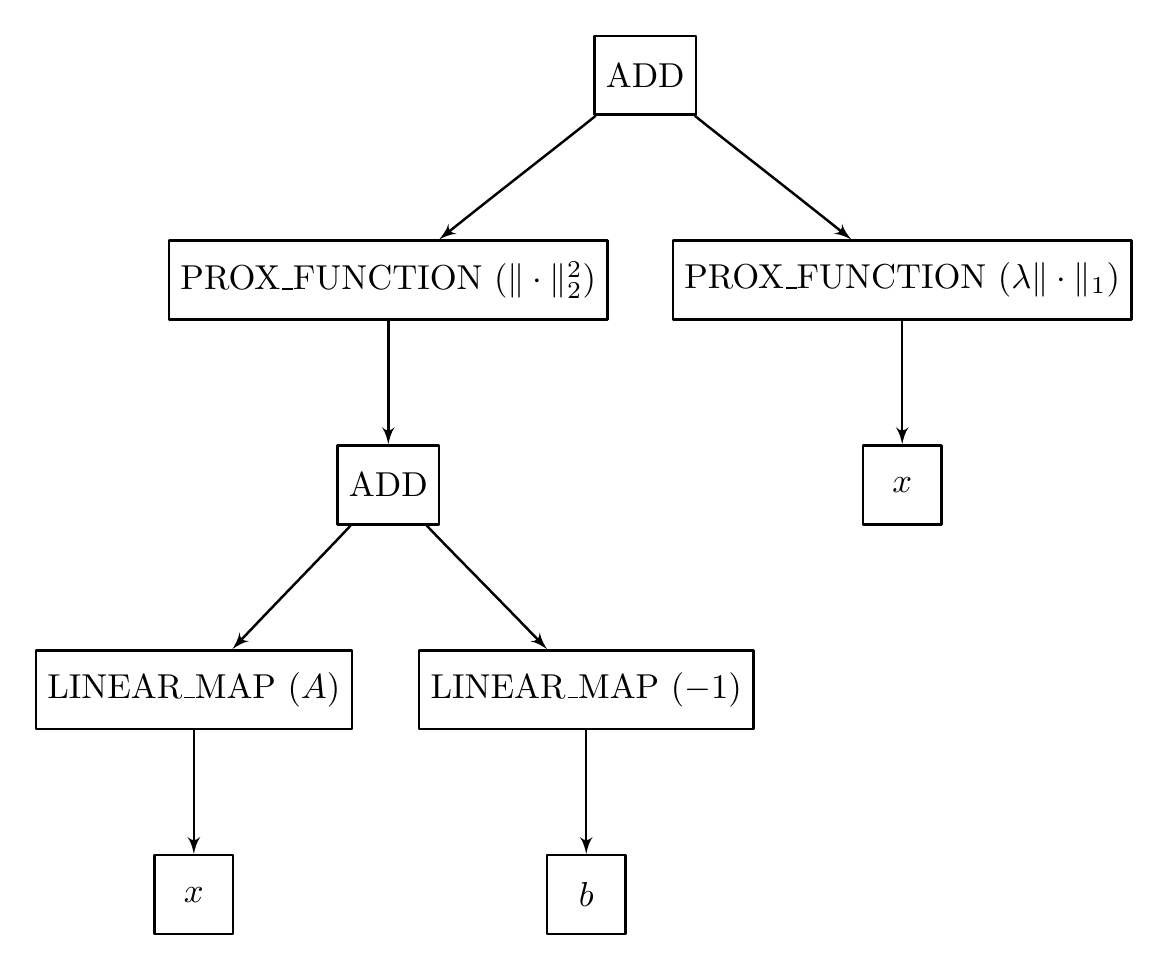}
  \caption{Original abstract syntax tree for lasso $\|Ax - b\|_2^2 +  \lambda\|x\|_1$
    (left) and converted to prox-affine form (right).}
  \label{fig-ast-lasso}
\end{figure}


The first stage of the Epsilon compiler transforms an arbitrary disciplined
convex problem into a prox-affine problem. Concretely, given an AST representing
the optimization problem in its original form (which may consist of any valid
composition of functions from the DCP library), this stage produces an AST with
a reduced set of nodes:
\begin{itemize}
\item \texttt{ADD}. The sum of its children $x_1 + \cdots + x_n$.
\item \texttt{PROX\_FUNCTION}. A prox-friendly function with proximal
  operator implementation in the Epsilon solver.
\item \texttt{LINEAR\_MAP}. A linear function with linear operator
  implementation in the Epsilon solver.
\item \texttt{VARIABLE}, \texttt{CONSTANT}. Each leaf of the AST is either a
  variable or constant.
\end{itemize}
As such, once the AST is transformed into prox-affine form, each node
corresponds to a function with a concrete proximal or linear operator
implementation from the library of operators described in Section
\ref{sec-atomic}; an example of this transformation is shown in Figure
\ref{fig-ast-lasso}.

\begin{algorithm}
  \caption{Conversion of AST to prox-affine form}
  \begin{algorithmic}[0]
    \Function{ConvertProx}{$rules$, $tree$}
    \Let{$rule$}{\Call{MatchRule}{$rules$, $tree$}}
    \If{$rule$ is nil}
    \State \Return{\Call{ConvertConic}{$tree$}}
    \EndIf

    \Let{$args$, $indicators$}{\Call{ConvertProxArguments}{$rule$, $tree$}}
    \Let{$output$}{\Call{CreateNode}{$rule$, $args$}}
    \For{$indicator$ in $indicators$}
    \Let{$output$}{\Call{Add}{$output$, \Call{ConvertProx}{$rules$, $indicator$}}}
    \EndFor
    \State \Return $output$
    \EndFunction
  \end{algorithmic}
  \label{alg-convert-prox}
\end{algorithm}


The transformation to prox-affine form is done in two passes over the AST
representing the optimization problem. In the
first pass, we convert all nodes representing linear functions to a nodes of
type \texttt{LINEAR\_MAP} which map directly to the linear operators available
in the Epsilon solver. In terms of ASTs representing linear functions, there are
two possibilities in any  DCP-valid input problem: 1) direct linear
transformaitons of inputs, such as the  unary \texttt{SUM} or variable argument
$\texttt{HSTACK}$ and 2) binary functions such as \texttt{MULTIPLY} which apply
a linear operator defined by a constant expression.  The former have
straightforward transformations to \texttt{LINEAR\_MAP} representations; for the
latter, DCP rules require that one of the arguments be constant and thus this
argument is evaluated and converted to a linear operator (typically, a sparse,
dense or diagonal matrix) and a \texttt{LINEAR\_MAP} node with single argument.

In the second pass, we complete the transformation to prox-affine form by
applying a set of prioritized rules, preferring to map ASTs onto a high-level
proximal operator implementation when available but falling back to conic
transformations, when necessary; the process is shown in Algorithm
\ref{alg-convert-prox}. In short, given an input tree (or subtree) along with a
set of rules for proximal operator transformations, we match the input against
these rules. If a matching proximal operator is found, we then transform the
function arguments (represented by subtrees of the original input tree) so that
they have valid form for composition with the proximal operator in question. In
doing so, the process may introduce auxiliary variables and additional indicator
functions; the behavior of \texttt{ConvertProxArguments} depends on the
requirements of the particular proximal operator (see Section
\ref{sec-proximal-operators}), but some common examples include:
\begin{itemize}
\item \emph{No-op.} The expression $\max\{-x, 0\}$ with variable $x$ is the
  hinge function  $f(x) = \max\{x, 0\}$, composed with the linear transformation
  $-I$. This represents a valid proximal operator so no further transformations
  are necessary and the argument $-x$ is returned as-is.
\item \emph{Epigraph transformation.} The expression $\|Ax - b\|_1$ with
  constants $A$, $b$ and variable $x$ matches a the proximal operator for the
  $\ell_1$-norm $\|\cdot\|_1$, but it cannot be composed with an arbitrary
  affine function given the set of proximal operators available in the Epsilon
  solver. Therefore, a new variable $y$ is introduced along with the
  constraint $I_0(Ax - b - y)$; $y$ is returned as argument, resulting in the new
  expression $\|y\|_1$.
\item \emph{Kronecker product splitting.} The expression $\|AXB - C\|_F^2$ with
  constants $A, B, C$ and variable $X$ matches the proximal operator for
  sum-of-squares, but evaluation would require factoring $F^TF
  + I$ where $F = B^T \otimes A$ which cannot be done efficiently given the
  linear operators available. Therefore, the compiler
  introduces a new variable $Z$, modifies the argument to be $\|AZ - C\|_F^2$
  and introduces the constraint $I_0(XB - Z)$.
\end{itemize}
Once the arguments have been transformed to the proper form, \texttt{CreateProx}
creates a \texttt{PROX\_FUNCTION} node (with attribute specifying \emph{which}
proximal operator implementation) and the transformed arguments as children. Any
indicators that were added by the argument conversion process are themselves
recursively converted to prox-affine form and the result is accumulated in the
output under an \texttt{ADD} node.

\subsection{Optimization and separation of prox-affine form}

Once the problem has been put in prox-affine form, the next stage of the
compiler transforms it to be separable in preparation for the solver.
As a given problem typically has many different separable prox-affine
representations, this stage must balance the
tradeoff between per-iteration computational complexity and the overall number
of iterations that will be required to solve a given problem. As an extreme
example, we can often split problems into a large number
of proximal operators which are very cheap to evaluate but will require a large
number of iterations. The theoretical analysis of convergence rates
for operator splitting algorithms is an active area of research (see
e.g. \cite{davis2014convergence,nishihara2015general,giselsson2015tight}) but
Epsilon follows the simple philosophy of minimizing the total number of
proximal operators in the separable form provided that each operator can be
evaluated efficiently. This is implemented with multiple passes on the
prox-affine form which split the problem as needed until it it satisfies the
constraints required for applying the operator splitting algorithm described in
Section \ref{sec-epsilon-solving}.

Although we previously described the prox-affine form generically, where we were
minimizing over a single variable $x$ and each term could potentially depend on
all variable $f_i(H_i(x))$, the reality is that for many problems the $x$
variables are already ``naturally'' partitioned to some extent (for instance,
this arises in epigraph transformations, where one function in the prox-affine
form will only depend on epigraph variables).  Thus, to be more concrete, we
introduce a partitioning of the variables $x = (x_1, \ldots, x_k)$ (where
here each $x_j$ is itself a vector of appropriate size, and let $\mathcal{J}_i$
denote the set of all variables that are used in the $i$th prox operator, i.e.
our optimization problem becomes
\begin{equation}
\minimize_{x_1,\ldots,x_k} \; \sum_{i=1}^n f_i(H_i(x_{\mathcal{J}_i})).
\end{equation}
The process of separation is effectively one of introducing ``copies'' of
variables until we reach a point that each objective term $f_i$ has a unique set
of variables, and the interactions between variables are captured entirely by
the explicit equality constraints.

Given the form above, we describe the optimization problem via a bipartite graph,
between nodes corresponding to objective functions $f_1, \ldots, f_N$ (plus
additional equality constraints, in the final form) and nodes corresponding to
variables $x_1,\ldots,x_k$.  An edge exists between $f_i$ and $x_j$ if $j \in
\mathcal{J}_i$, i.e. if the function uses that variable.  By applying a
sequence of transformations, we will introduce new variables and new equality
constraints that will put the problem into a separable prox-affine form.

\paragraph{Definition of equivalence transformations.}

Specifically, the compiler sequentially executes a series of transformations to
put the problem in separable prox-affine form:
\begin{enumerate}
\item \emph{Move equality indicators.} The first compiler stage (``Conversion to
  prox-affine form'', see Section \ref{sec-compiler-conversion}) produces a
  single expression for the objective which includes all constraints via
  indicators; due to the nature of the transformations, many equality
  constraints are ``simple'' (e.g. involving $I$ or $-I$) and can thus be moved
  to actual constraints in the separable form. This pass performs these
  modifications based on the linear map associated with the edges
  corresponding to variables in each equality constraint, splitting expressions
  when necessary. For example, an objective term $I_0(Ax + y + z)$ is
  transformed to a new objective term $I_0(Ax - w)$ and the constraint $w + y +
  z = 0$.
\item \emph{Combine objective terms.} The basic properties of proximal operators
  (see e.g. \cite{parikh2013proximal}) allow simple functions like $c^Tx$ and $\|x - b\|_2^2$ to be
  combined with other terms, reducing the number of proximal operators needed in
  the separable prox-affine form. This pass combines these terms assuming there
  is another objective term which includes the same variable.
\item \emph{Add variable copies and consensus constraints.} The final pass
  guarantees that the objective is separable by introducing variable copies and
  consensus constraints. For example, the objective $f(x) + g(x)$ is
  transformed to $f(x_1) + g(x_2)$ and the constraint $x_1 = x_2$ is added.
\end{enumerate}

\begin{figure}
  \centering
  \includegraphics[scale=0.8]{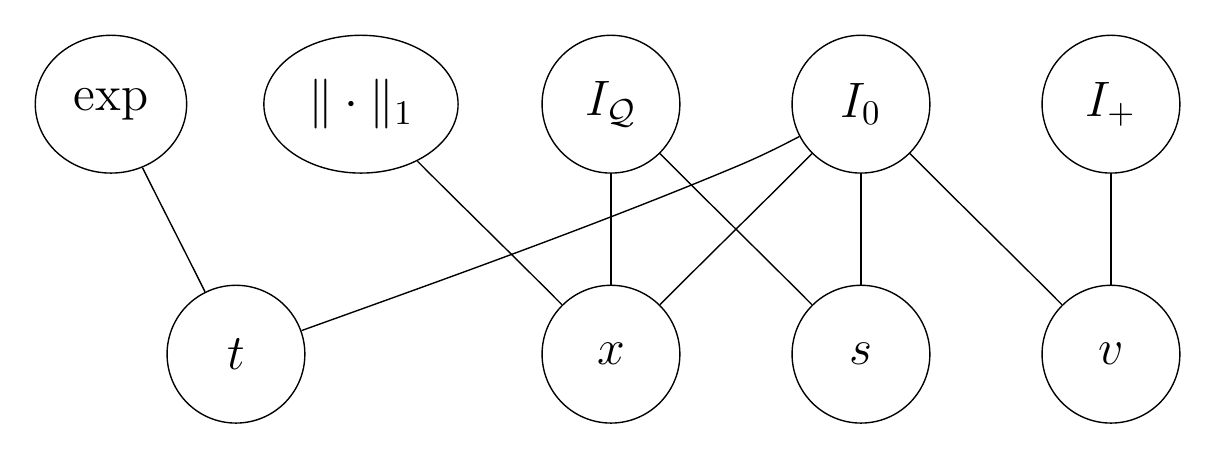} \\
  \includegraphics[scale=0.8]{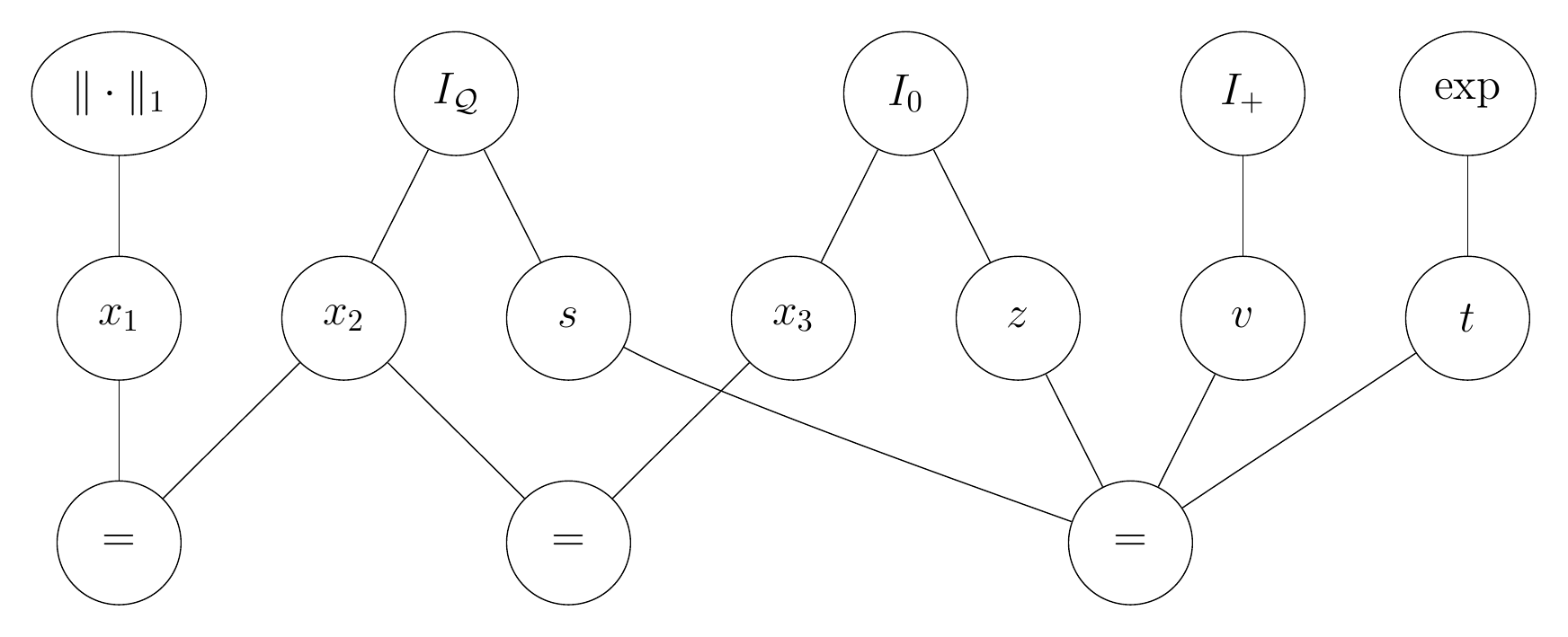}
  \caption{Compiler representation of prox-affine form as bipartite graph for
    $\exp(\|x\|_2 + c^Tx) + \|x\|_1$ after conversion pass (top, see equation
    \ref{eq-exp-prox-affine}) and after separation pass (bottom, see equation
    \ref{eq-exp-sep-prox-affine}).}
  \label{fig-graph-exp}
\end{figure}


For illustration purposes, consider the problem
\begin{equation}
  \minimize_x \;\; \exp(\|x\|_2 + c^Tx) + \|x\|_1.
\end{equation}
The first compiler stage converts this problem to prox-affine form by introducing
auxiliary variables $t,s,v$, along with three cone constraints and two
prox-friendly functions:
\begin{equation}
  \label{eq-exp-prox-affine}
  \minimize_{x,t,s,v} \;\; \exp(t) + \|x\|_1 + I_{\mathcal{Q}}(x,s) + I_0(s +
  c^Tx - t - v) + I_+(v),
\end{equation}
where $I_\mathcal{Q}$ denotes the indicator of the second-order cone, $I_0$
the zero cone and $I_+$ the nonnegative orthant. The problem in this form is the
input for the second stage which constructs the bipartite graph shown in Figure
\ref{fig-graph-exp} (top). The problem is then transformed to have separable
objective with many of the terms in $I_0$ (those with simple linear maps) move
to the constraint
\begin{equation}
  \label{eq-exp-sep-prox-affine}
  \begin{split}
  \minimize_{x_1,x_2,x_3,t,s,v} \;\; & \exp(t) + \|x_1\|_1 +
  I_{\mathcal{Q}}(x_2,s) + I_0(c^Tx_3 - z) + I_+(v) \\
  \subjectto \;\; & s + z - t - v = 0 \\
  & x_1 = x_2 \\
  & x_2 = x_3
  \end{split}
\end{equation}
and variables $x_1, x_2, x_3$ are introduced along with consensus
constraints. The bipartite graph for the final output from the compiler, a
problem in separable prox-affine form, is shown in Figure \ref{fig-graph-exp}
(bottom).

\subsection{Solving problems in prox-affine form}
\label{sec-epsilon-solving}

Once the problem has been put in separable prox-affine form, the Epsilon solver applies the
ADMM-based operator splitting algorithm using the library of proximal and linear
operators. The implementation details of each operator are abstracted from the
high-level algorithm which applies the operators only through a common
interface providing the basic mathematical operations required. Next we give
a mathematical description of the operator splitting algorithm itself while the
computational details of individual proximal and linear operators are discussed
in Section \ref{sec-atomic}.

The Epsilon solver employs a variant of ADMM to solve problems in the separable
prox-affine form
\begin{equation}
  \begin{split}
    \minimize_{x_i,\ldots,x_N} \;\; & \sum_{i=1}^N f_i(H_i(x_i)) \\
    \subjectto \;\; & \sum_{i=1}^N A_i(x_i) = b.
  \end{split}
\end{equation}
This approach can be motivated by considering the augmented Lagrangian
\begin{equation}
  L_\lambda(x_1,\ldots,x_N,y) = \sum_{i=1}^N f_i(H(x_i)) + y^T(Ax - b) +
  (1/2\lambda) \| Ax - b \|_2^2
\end{equation}
where $y$ is the dual variable, $\lambda \ge 0$ is the augmented Lagrangian
penalization parameter, and $Ax = \sum_{i=1}^NA_i(x_i)$. The ADMM method applied
here results in the Gauss-Seidel updates\footnote{specifically, the update for
$x_i^{(k+1)}$ depends on $x_j^{(k+1)}$ for $j<i$ and $x_j^{(k)}$ for $j>i$} with
\begin{equation}
  \begin{split}
    x_i^{k+1} &:=  \argmin_{x_i} \lambda f_i(H_i(x_i)) + \frac{1}{2} \left \|
    \sum_{j < i}A_j(x_j^{k+1}) + A_i(x_i) + \sum_{j > i} A_j(x_j^{k}) - b + u^k \right\|_2^2 \\
    u^{k+1}   &:= u^k + Ax^{k+1} - b
  \end{split}
\end{equation}
where we have $u = \lambda y$ is the scaled dual variable. Critically, the
$x_i$-updates are applied using the (generalized) proximal operator,
let
\begin{equation}
v_i^k = b - u^k - \sum_{j < i}A_j(x_j^{k+1}) - \sum_{j > i}
A_j(x_j^k)
\end{equation}
and we have
\begin{equation}
  \label{ir-xi-update}
    x_i^{k+1} :=  \argmin_{x_i} \lambda f_i(H_i(x_i)) + (1/2)\|A(x_i) - v_i^k
    \|_2^2 = \prox_{\lambda f_i \circ H_i, A_i}(v_i^k)
\end{equation}
The ability of the solver to evaluate the
generalized proximal operator efficiently will depend on $f_i$ and $A_i$ (in the
most common case $A_i^TA_i = \alpha I$, a scalar matrix, which can be handled by
any proximal operator); it is the responsibility of the compiler to ensure that
the prox-affine problem has been put in the required form such that these
evaluations map to efficient implementations from the proximal operator library.

\section{Fast atomic operators}
\label{sec-atomic}


The Epsilon solver contains a library of fast linear and
proximal operators which implement the operations required for a high-level
iterative algorithm such as the ADMM variant described in the previous
section. In particular, the proximal operator library directly implements a large number of
the functions available to disciplined convex programming frameworks reducing
the need for extensive transformation before solving an optimization
problem. As the evaluation of proximal operators as well as the
operations required by high-level algorithms rely heavily on linear
operators, Epsilon also provides a library of efficient linear operators (and a
system for composing them), extending beyond the standard dense/sparse matrices
typically found in generic convex solvers.

\subsection{Linear operators}
\label{sec-linear-operators}

In general, the computation required for solving convex optimization often typically
depends heavily on the application of linear operators. Most commonly these
linear operators are implemented with sparse or dense matrices which explicitly
represent the coefficients of the linear transformation. Clearly, in many cases
this can be inefficient (see e.g. \cite{diamond2015convex} and the references therein) and as
such, we abstract the notion of a linear operator allowing for other
implementations which can often be far more efficient than direct matrix
representation.

A motivating example that arises in many applications is the use of matrix-valued variables. As
intermediate representations for convex programming (both prox-affine and conic
forms) typically reduce optimization problems over matrices to ones over
vectors, matrix products naturally give rise to the Kronecker product. For
example, consider the expression $AX$ where $A \in \mathbb{R}^{m \times n}$ is a
dense constant and $X
\in \mathbb{R}^{n \times k}$ is the optimization variable; we
vectorize this product with $\vect(AX) = (I \otimes A)\vect(X)$ where $\otimes$ denotes the
Kronecker product. Representing the Kronecker product as a sparse matrix is not
only space inefficient (i.e. requiring $A$ to be repeated $k$
times) but can also be extremely costly to factor. In particular, the naive
approach requires a sparse factorization of a $km \times kn$ matrix as opposed
to factoring a dense $m \times n$ matrix $A$ directly. Explicitly
maintaining the Kronecker product structure provides a mechanism
for avoiding this unnecessary computational cost.

Epsilon augments the standard
sparse/dense matrices with dedicated implementations for diagonal matrices, scalar
matrices, the Kronecker product as well as composite types representing a sum or
product:
\begin{itemize}
\item \emph{Dense matrix}. A dense matrix $A \in \mathbb{R}^{m \times n}$ with
  $O(mn)$ storage.
\item \emph{Sparse matrix}. A sparse matrix $A \in \mathbb{R}^{m \times n}$ with
  O(\# nonzeros) storage.
\item \emph{Diagonal matrix}. A diagonal matrix $A \in \mathbb{R}^{n \times n}$
  with $O(n)$ storage.
\item \emph{Scalar matrix}. A scalar matrix $\alpha I \in \mathbb{R}^{n \times
  n}$ with  $\alpha \in \mathbb{R}$ and $O(1)$ storage.
\item \emph{Kronecker product}. The Kronecker product of $A \in \mathbb{R}^{m
  \times n}$ and $B \in \mathbb{R}^{p \times q}$, representing a linear map
  $\mathbb{R}^{nq} \to \mathbb{R}^{mp}$ where $A$ and $B$ can themselves be any
  linear operator type.
\item \emph{Sum}. The sum of linear operators $A_1 + \cdots + A_N$.
\item \emph{Product}. The product of linear operators $A_1\cdots A_N$.
\item \emph{Abstract}. An abstract linear operator which cannot be combined with
  any of the basic types. This is used to represent factorizations, see below.
\end{itemize}

Each linear operator type $A$ supports the following operations:
\begin{itemize}
\item \emph{Apply}. Given a vector $x$, return $y = Ax$.
\item \emph{Transpose}. Return the linear operator $A^T$.
\item \emph{Inverse}. Return the linear operator $A^{-1}$.
\end{itemize}
The transpose operation returns a linear operator of the same type whereas the
inverse operation may return a linear operator of a different type. The inverse
operation is undefined for non-invertible linear operators and in practice is
intended to be used in contexts where the linear transformation is known to be
invertible.

In addition, linear operators also support binary operations for
sum $A + B$ and product $AB$. This requires a system for type conversion,
the basic rules for which are described with an ordering of the types
corresponding to their sparsity ($Dense > Sparse > Diagonal > Scalar$); in order
to combine any two of these types, we first promote the sparser type to
to the denser type. For Kronecker products, type
conversion depends on the arguments: the sum of two Kronecker products can be
combined if one of the arguments is equivalent, e.g.
\begin{equation}
  A \otimes B + A \otimes C = A \otimes (B + C)
\end{equation}
while the product can be combined if the arguments have matching dimensions,
i.e.
\begin{equation}
  (A \otimes B)(C \otimes D) = AC \otimes BD
\end{equation}
if we can form $AC$ and $BD$. In either case, if a combination is possible it
will be performed along with the appropriate type conversion for the
sum/product of the arguments themselves. If a combination is not possible
according to these rules, the resulting type will instead be a $Sum$
or $Product$ composed of the arguments.

\subsection{Proximal operators}
\label{sec-proximal-operators}


The second class of operators that
form the basis for Epsilon are the proximal operators, described above in
Section \ref{sec-background-prox}. As seen in the resulting description of the
ADMM algorithm, each individual solution over the $x_i$ variables can be
represented via an operator
\begin{equation}
x^{k+1}_i := \argmin_{x_i} \lambda f(H(x_i)) + \frac{1}{2}\|A_i(x_i) - v_i^k\|^2_2
\end{equation}
for some value of $v_i^k$ (which naturally depends generally on the dual variables
for the equality constraints involving $x_i$ as well as the other $x_j$
variables). This is exactly the generalized proximal operator $x^{k+1}_i =
\prox_{\lambda f \circ H_i,A_i} (v)$. Indeed, the Epsilon compiler uses precisely the set of available fast
proximal operators to reduce the convex optimization problems to fast forms
relative to the corresponding cone problem; while any of the problems can be solved by
reducing everything to these conic form (and thus using only proximal
operators corresponding to cone projections), the speed of the solver crucially
depends on the ability to evaluate a much wider range of these proximal operators
efficiently.


It is well-known that many proximal
operators have closed form solutions that can be solved much more quickly than
general optimization problems (see e.g. \cite{parikh2013proximal} for a review).
As the main goal of this current paper is not to develop new proximal operators,
we merely highlight several of the operators included in Epsilon along with a
general description of the methods used to solve them.  The proximal operators
group roughly into three categories: elementwise, vector, and matrix operators,
for functions $f$ of scalar inputs, vector inputs, and matrix inputs
respectively.  Note that this is not a perfect division, because several vector
matrix functions are simply sums of corresponding scalar functions, e.g.,
$\|x\|_2^2 = \sum_{i=1}^n x_i^2$, but instead we use vector or matrix designations to
refer to functions that \emph{cannot} be decomposed over their inputs.

\begin{itemize}

\item \textbf{\textsf{Exact linear, quadratic, or cubic
equations.}}  Several prox operators can be solved using exact solutions to
their gradient conditions.  For instance, the prox operator of $f(x) = (ax -
b)^2$ is given by the solution to the gradient equation $2a(ax - b) + x = 0$,
which is a simple linear equation.  Similarly, the prox operators for $f(x) =
-\log x$ and $f(x) = 1/x$ have gradient conditions that are given by quadratic
and cubic equations respectively.  Here some care must be taken to ensure that
we select the correct of two or three possible solutions, but this can be
ensured analytically in many cases or simply by checking all solutions in the
worst case.  For vector functions, the proximal operator of a least-squares
objective is also an instance of this case, though here of course the complexity
requires a matrix inversion rather than a scalar linear equation.

\item \textbf{\textsf{Soft thresholding.}}  The absolute value and related
functions can be solved via the soft thresholding operator.  For example,
the proximal operator of $f(x) = |x|$ is given by
\begin{equation}
\prox_{\lambda f}(v) = \left \{ \begin{array}{ll} v - \lambda & \mbox{ if } v
> \lambda \\ 0 & \mbox{ if } -\lambda \leq v \leq \lambda \\ v + \lambda & \mbox
{ if } v < -\lambda. \end{array} \right .
\end{equation}

\item \textbf{\textsf{Newton's method.}}  Several proximal operators for smooth
functions have no easily computed closed form solution.  Nonetheless, in this
case Newton's method can be used to find a solution in reasonable time.  For
elementwise functions, for instance, these operations are relatively efficient,
because in practice a very small number of Newton iterations are needed to reach
numerical precision.  For example, the proximal operator of the logistic
function $f (x) = \log (1 + \exp (x))$ has no closed form solution, but can
easily computed with Newton's method.

\item \textbf{\textsf{Projection approaches.}} The proximal operator for several
functions can be related to the projection on to a set.  In the simplest
case, the proximal operator for an indicator of a set is simply equal to the
projection onto the set, giving prox operators for cone projections.
However, additional proximal methods derive from Moreau decomposition \cite{moreau1962fonctions},
which states that
\begin{equation}
v = \prox_f(v) + \prox_{f^*}(v)
\end{equation}
where $f^*$ denotes the Fenchel conjugate of $f$.  For example, the proximal
operator for the $\ell_\infty$-norm, $f(x) = \|x\|_\infty$ is given by
\begin{equation}
\prox_{f}(v) = v - \mathcal{P}_{\|x\|_1 \leq 1}(v)
\end{equation}
using the relation that the Fenchel conjugate of a normal is equal to the
indicator of of the dual norm ball.  The projection on to the $\ell_1$-ball can
be accomplished in $O(n)$ time using a median of medians algorithm \cite{duchi2008efficient}.

\item \textbf{\textsf{Orthogonally invariant matrix functions.}} If a matrix
input function is defined solely by the singular values (or eigenvalues) of a
matrix, then the proximal operator can be computed using the singular value
decomposition (or eigenvalue decomposition) of that matrix.  Typically, the
running time of these proximal operators are dominated by the cost of the
decomposition itself, making them very efficient for reasonably-sized matrices.
For example, the $\log\det$ function can be written as
\begin{equation}
\log \det(X) = \sum_{i=1}^n \log \lambda_i
\end{equation}
so its proximal operator can be given by
\begin{equation}
\prox_{\log\det}(X) = U \prox_{\log}(\Lambda) U^T
\end{equation}
where $\prox_{\log}(\Lambda)$ is shorthand for applying the proximal operator
for the negative log function to each diagonal element of the eigenvalues
$\Lambda$. The prox operator for the $\log$ function can itself be solved via a
quadratic equation, so computing the inner prox term is only an $O(n)$
operation and the runtime is dominated by the $O(n^3)$ cost of computing the
eigenvalue decomposition.

\item \textbf{\textsf{Special purpose algorithms.}}  Finally, though we cannot
enumerate these broadly, several proximal operators have special purpose fast
solvers for these particular types of operators.  A particularly relevant
example is the fused lasso $f(x) = \|D x\|_1$ where $D$ is the first order
difference operator; this proximal operator can be solved efficiently via an $O
(n)$ dynamic programming approach \cite{johnson2013dynamic}.

\end{itemize}


\begin{landscape}
\begin{minipage}{9in}

\begin{center}
\begin{footnotesize}
\begin{tabular}{|c|c|c|c|c|}
\hline
\multicolumn{3}{|c|}{\textbf{Function}} &
\multicolumn{2}{|c|}{\textbf{Proximal operator}}
\\ \hline
\textbf{Category} & \textbf{Atom} & \textbf{Definition} & \textbf{Method} &
\textbf{Complexity} \\
\hline

\multirow{4}{*}{Cone}
& Zero & $f(x) = I_0(Ax - b)$, $A \in \mathbb{R}^{m \times n}$
& subspace projection & $O(mn^2 + n^3)$ \\

& Nonnegative orthant & $f(x) = I_+(x)$
& positive thresholding & $O(n)$ \\

& Second-order cone & $f(x,t) = I_{\mathcal{Q}}(x,t)$, $x \in \mathbb{R}^n$, $t \in \mathbb{R}$
& projection & $O(n)$ \\

& Semidefinite cone & $f(X) = I_\succeq(X), X \in \mathbb{S}^n$
& positive thresholding on $\lambda(X)$ & $O(n^3)$ \\

\hline
\multirow{12}{*}{Elementwise $x,y \in \mathbb{R}$ }
& Absolute & $f(x) = |x|$
& soft thresholding & $O(n)$ \\

& Square & $f(x) = x^2$ &
linear equation & $O(n)$ \\

& Hinge & $f(x) = \max\{x,0\}$
& soft thresholding & $O(n)$ \\

& Deadzone & $f(x) = \max\{|x| - \epsilon, 0\}, \; \epsilon \geq 0$
& soft thresholding & $O(n)$ \\

& Quantile & $f(x) = \max\{\alpha x,(\alpha-1)x\}, \; 0 \leq \alpha \leq 1$
& asymmetric soft thresholding & $O(n)$ \\

& Logistic & $f(x) = \log(1+\exp(x))$
& Newton & $O(n)\cdot(\mbox{\# Newton})$ \\

& Inverse positive & $f(x) = 1/x, \; x\geq 0$
& Newton & $O(n)\cdot(\mbox{\# Newton})$ \\

& Negative log & $f(x) = -\log(x), \; x\geq 0$
& quadratic equation & $O(n)$ \\

& Exponential & $f(x) = \exp(x)$
& Newton & $O(n)\cdot(\mbox{\# Newton})$ \\

& Negative entropy & $f(x) = x \cdot \log(x), \; x \geq 0$
& Newton & $O(n)\cdot(\mbox{\# Newton})$ \\

& KL Divergence & $f(x,y) = x \cdot \log(x/y), \; x,y \geq 0$
& Newton & $O(n)\cdot(\mbox{\# Newton})$ \\

& Quadratic over linear & $f(x,y) = x^2/y, \; y \geq 0$
& cubic equation & $O(n)$ \\

\hline
\multirow{6}{*}{Vector $x \in \mathbb{R}^n$}

& $\ell_1$-norm & $f(x) = \|x\|_1$
& soft thresholding & $O(n)$  \\

& Sum-of-squares & $f(x) = \|Ax - b\|^2_2$, $A \in \mathbb{R}^{m \times n}$
& normal equations & $O(mn^2 + n^3)$\footnote[$\dagger$]{Or $O(nm^2 + m^3)$ if
$m < n$, and often lower if $A$ is sparse.  Furthermore, the cost can be
amortized over multiple evaluations for the same $A$ matrix,
we can compute a Cholesky factorization once in this time, and solve
subsequent iterations in $O(mn + n^2)$ time.} \\

& $\ell_2$-norm & $f(x) = \|x\|_2$
& group soft thresholding & $O(n)$ \\

& $\ell_\infty$-norm & $f(x) = \|x\|_\infty$
& median of medians & $O(n)$ \\

& Log-sum-exp & $f(x) = \log \left(\sum_{i=1}^n \exp(x_i)\right)$
& Newton & $O(n)\cdot(\mbox{\# Newton})$ \\

& Fused lasso & $f(x) = \sum_{i=1}^{n-1} |x_i - x_{i+1}|$ & dynamic programming \cite{johnson2013dynamic} & $O(n)$ \\

\hline
\multirow{3}{*}{Matrix}

& Negative log det & $f(X) = -\log \det(X), X \in \mathbb{S}^n$
& quadratic equation on $\lambda(X)$ & $O(n^3)$ \\

& Nuclear norm & $f(X) = \|\sigma(X)\|_1 , X \in \mathbb{R}^{m \times n}$
& soft thresholding on $\sigma(X)$ & $O(n^3)$ \\

& Spectral norm  & $f(X) = \|\sigma(X)\|_\infty , X \in \mathbb{R}^{m \times n}$
& median of medians on $\sigma(X)$ & $O(n^3)$ \\

\hline
\end{tabular}
\end{footnotesize}
\end{center}
\end{minipage}
\end{landscape}

\section{Examples and numerical results}
\label{sec-examples}

In this section we present several examples of convex problems from statistical
machine learning and compare Epsilon to existing approaches on a library of
examples from several domains. At present, Python integration is provided
(Matlab, R, Julia versions are planned) with CVXPY \cite{diamond2015cvxpy}
providing the environment for constructing the disciplined convex programs and
performing DCP validation. Epsilon effectively serves as a solver for
CVXPY although behind the scenes the CVXPY/Epsilon interface is somewhat
different than for other solvers as Epsilon compiler implements its own
transformations on the AST. Nonetheless, from a user perspective problems are
specified in the same syntax, a high-level domain specific language of which we
give several examples in this section. Epsilon is open source and available at
\url{http://github.com/mwytock/epsilon/}, including the code for all examples and
benchmarks discussed here.

As Epsilon integrates with CVXPY, we make the natural comparison between Epsilon
and the existing solvers for CVXPY which use the conic form. In particular, we
compareEpsilon to ECOS \cite{domahidi2013ecos}, an interior point method and SCS
\cite{o2013operator}, the ``splitting conic solver''. In general, interior point methods achieve
highly accurate solutions but have trouble scaling to larger problems and so
it is unsurprising that Epsilon is able to solve problems to moderate
accuracy several orders of magnitude faster than ECOS (this is also the case
when comparing ECOS and SCS, see \cite{o2013operator}). On the other hand, SCS
employs an operator splitting method that is similar in spirit to the Epsilon
solver, both being variants of ADMM; the main difference between the two is in the
intermediate representation and set of available proximal operators: SCS solves the conic
form produced by CVXPY with cone and subspace projections (using sparse
matrices) while Epsilon solves the prox-affine form using the
much larger library of fast proximal and linear operators described in Section
\ref{sec-atomic}. In practice, this has significant impact with Epsilon
achieving the same level of accuracy as SCS an order of magnitude faster (or
more) on many problems.

In what follows we examine several examples in detail followed by results on a
library of 20+ problems common to applications in statistical machine learning
and other domains. In the detailed examples, we start with a
complete specification of the input problem (a few lines of Python in the CVXPY
grammar), discuss the transformation by the Epsilon compiler to an abstract syntax tree
representing a prox-affine problem and explore how the Epsilon solver scales
relative to conic solvers. When printing the AST for individual problems, we
adopt a concise functional form which is a serialized version of the abstract
syntax trees shown in Figure \ref{fig-ast-lasso} and is the internal
representation of the Epsilon compiler, as well as input to the Epsilon solver.

\subsection{Lasso}

We start with the lasso problem
\begin{equation}
\minimize_\theta \;\; (1/2)\|X\theta - y\|_2^2 + \lambda \|\theta\|_1,
\end{equation}
where $X \in \mathbb{R}^{m \times n}$ contains input features, $y \in
\mathbb{R}^m$ the outputs, and $\theta \in \mathbb{R}^n$ are the model
parameters. The regularization parameter $\lambda \ge 0$ controls the tradeoff between data
fit and sparsity in the solution---the lasso is especially useful in the
high-dimensional case where $m \le n$ as sparsity effectively controls the
number of free parameters in the model, see \cite{tibshirani1996regression} for
details.

In CVXPY, the lasso can be written as
\begin{verbatim}
theta = Variable(n)
f = sum_squares(X*theta - y) + lam*norm1(theta)
prob = Problem(Minimize(f))
\end{verbatim}
where \texttt{sum\_squares()}/\texttt{norm1()} functions correspond directly to
the two objective terms. In essence this problem is already in prox-affine form
with proximal operators for $\|X\theta - y\|_2^2$ and $\|\theta\|_1$; thus the
prox-affine AST produced by the Epsilon compiler merely adds an additional
variable and equality constraint to make the objective separable:
\begin{verbatim}
objective:
  add(
    sum_squares(add(const(a), scalar(-1.00)*dense(B)*var(x))),
    norm1(var(y)))

constraints:
  zero(add(var(y), scalar(-1.00)*var(x)))
\end{verbatim}
Note that in the automatically generated serialization of the AST, variable/constant
names are auto-generated and do not necessarily correspond to user input. In
this particular example \texttt{a}, \texttt{B} correspond to the constants $y$, $X$
from the original problem while \texttt{x}, \texttt{y} correspond to two copies
of the optimization variable $\theta$.

\begin{figure}
  \centering
  \includegraphics[scale=0.5]{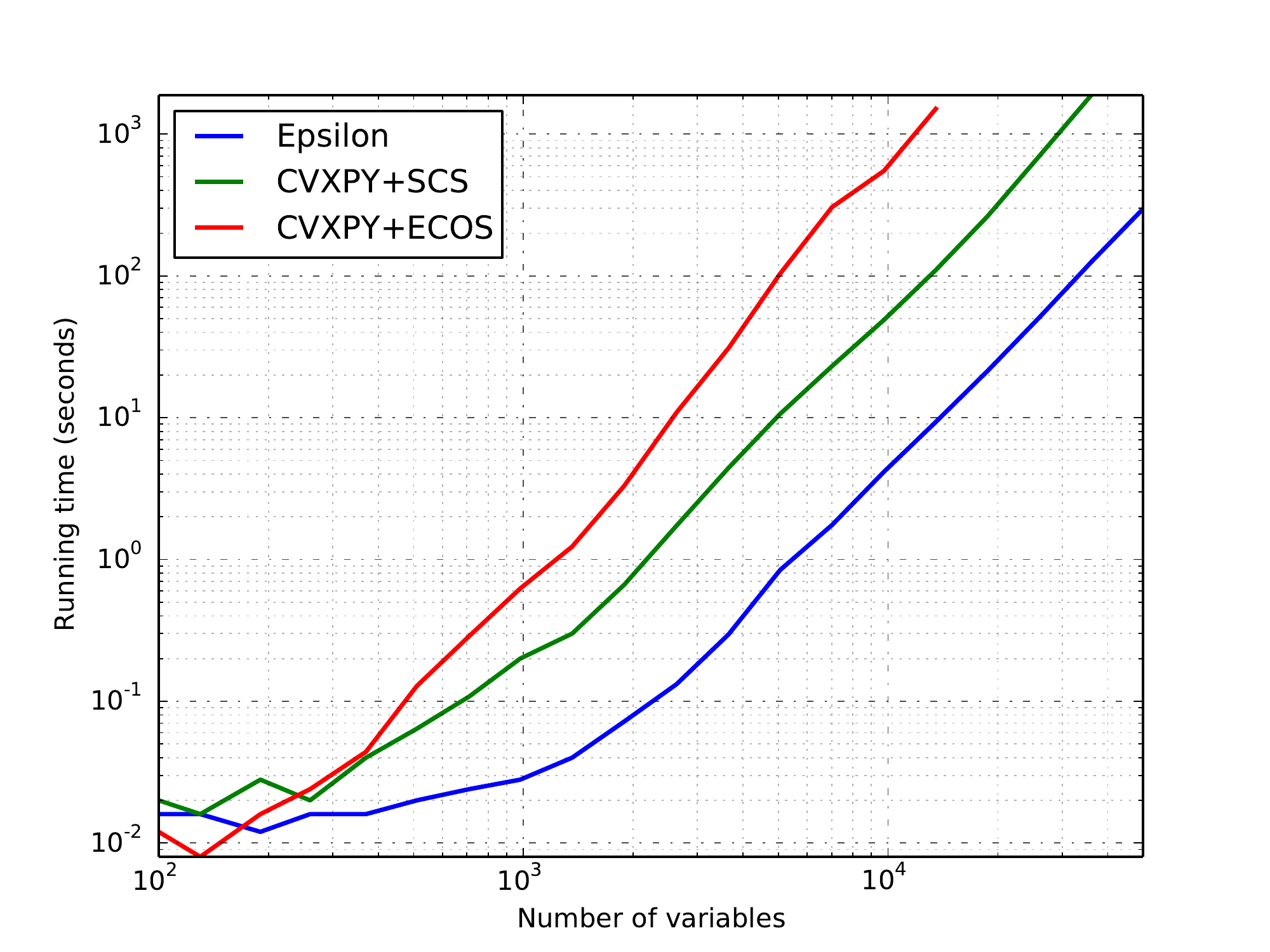}
  \caption{Comparison of running times on lasso example with dense $X \in
    \mathbb{R}^{m \times 10m}$.}
  \label{fig-scale-lasso}
\end{figure}

Computationally, it is the evaluation of
the sum-of-squares proximal operator which dominates the running time
for Epsilon as it requires solving the normal equations. However, the cost of
the factorization required is amortized as this step can be performed once
before the first iteration of the algorithm, as discussed in Section \ref{sec-proximal-operators}. In
Figure \ref{fig-scale-lasso}, we compare the running time of Epsilon to
CVXPY+SCS/ECOS on a sequence of problems with dense $X \in \mathbb{R}^{m \times
  10m}$. Epsilon (representing $X$ as a dense linear operator) performs a dense
factorization while SCS embeds $X$ in a large sparse matrix and performs a
sparse factorization (as it does with all problems). The difference in these
factorizations explains the difference in running time as the time spent
performing the actual iterations is negligible for both methods.

\subsection{Multivariate lasso}

In this example, we apply lasso to the multivariate regression setting where the
output variable is now vector as opposed to a scalar in univariate regression. In
particular,
\begin{equation}
  \minimize_\Theta \;\; (1/2)\|X\Theta - Y\|_F^2 + \lambda \|\Theta\|_1
\end{equation}
where $X \in \mathbb{R}^{m \times n}$ are input features, $Y \in \mathbb{R}^{m
  \times k}$ represent the $k$-dimensional response variable and the
optimization variable is now a matrix $\Theta \in \mathbb{R}^{n \times k}$,
representing the parameters of the multivariate regression model. The Frobenius norm
$\|\cdot\|_F$ is the $\ell_2$-norm applied elementwise to a matrix and
here $\|\cdot\|_1$ is also interpreted elementwise.

The CVXPY problem specification for the multivariate lasso is virtually
identical to the standard lasso example
\begin{verbatim}
Theta = Variable(n,k)
f = sum_squares(X*Theta - Y) + lam*norm1(Theta)
prob = Problem(Minimize(f))
\end{verbatim}
with the only change being the replacement of vectors \texttt{y}, \texttt{theta}
with their matrix counterparts \texttt{Y}, \texttt{Theta} (by convention, we
denote matrix-valued variables with capital letters). As a result, when the
Epsilon compiler transforms this problem to the prox-affine AST,
\begin{verbatim}
objective:
  add(
    sum_squares(add(kron(scalar(1.00), dense(A))*var(X), scalar(-1.00)*const(B))),
    norm1(var(Y)))

constraints:
  zero(add(var(Y), scalar(-1.00)*var(X)))
\end{verbatim}
the matrix-valued optimization variable results in the \texttt{kron} linear
operator appearing as argument to the \texttt{sum\_squares} proximal
operator. This corresponds to the specialized Kronecker product linear operator
implementation with (in this case) $O(mn)$ for the dense data matrix $X$.

\begin{figure}
  \centering
  \includegraphics[scale=0.5]{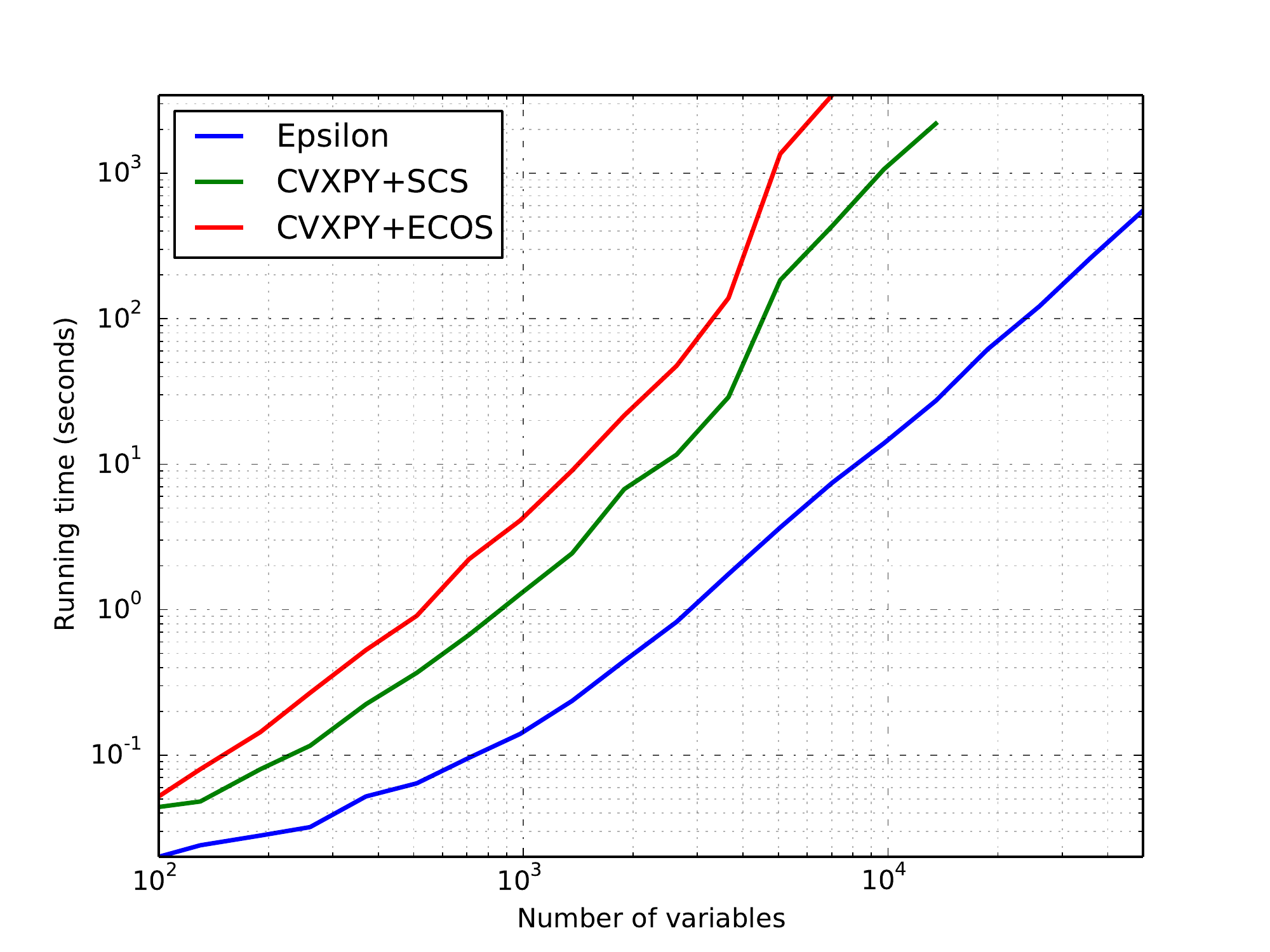}
  \caption{Comparison of running times on the multivariate Lasso example with dense $X \in
    \mathbb{R}^{m \times 10m}$ and $k = 10$.}
  \label{fig-scale-mv-lasso}
\end{figure}

Although it is a simple to change to the problem specification to
apply Lasso in the multivariate case, the new problem results in a substantially
different computational running
times for the solvers considered. In Figure \ref{fig-scale-mv-lasso} we show
the running times of each approach on a sequence of problems with $X \in
\mathbb{R}^{m \times 10m}$ and $Y \in \mathbb{R}^{m \times 10}$; where as on
standard Lasso Epsilon was roughly 10x faster than SCS, now the gap is
closer to 100x, e.g. for a problem with $1.35 \times 10^4$ variables, SCS
requires 2192 seconds vs. 27 seconds for Epsilon. The reason for this difference
is due to the representation of the linear operator required for solving the
normal equations for the least squares term
\begin{equation}
  (1/2) \|(I_k \otimes X)\vect(\Theta) - \vect(Y) \|_2^2.
\end{equation}
Since SCS is restricted to representing linear operators as
sparse matrices, it must instantiate the Kronecker product explicitly
(replicating the $X$ matrix $10$ times) and factor the resulting matrix with
sparse methods. In contrast, Epsilon represents the Kronecker product directly
(using the \texttt{kron} linear operator) and applies dense factorization
methods without any unnecessary replication.

\subsection{Total variation}

In the previous lasso examples, we employ the
$\ell_1$-norm to estimate a sparse set of regression coefficients---a natural
extension to this idea is to incorporate a notion of structured sparsity. The
fused lasso problem (originally proposed in \cite{tibshirani2005sparsity})
\begin{equation}
  \minimize_\theta \;\; (1/2) \|X\theta - y\|_2^2 + \lambda_1 \|\theta\|_1 +
  \lambda_2 \sum_{i=1}^{n-1} \|\theta_{i+1} - \theta_i\|_1
\end{equation}
employs total variation regularization (originally proposed in
\cite{rudin1992nonlinear}) to encourage the \emph{differences} of the
coefficient vector $\theta_{i+1} - \theta_i$ to be sparse. Such structure
naturally arises in problems where the dimensions of the coefficient vector
correspond to vertices in a chain or grid, see
e.g. \cite{xin2014efficient,adhikari2015high} for example applications.

For total variation problems, CVXPY
provides a function \texttt{tv()} which makes the problem specification
concise:
\begin{verbatim}
theta = Variable(n)
f = sum_squares(X*theta - y) + lam1*norm1(theta) + lam2*tv(theta)
prob = Problem(Minimize(f)).
\end{verbatim}
In conic form this penalty is transformed to a set of a linear constraints which
involve the first order differencing operator (to be defined shortly); however,
Epsilon includes a direct proximal operator implementation of the
total variation penalty and thus the compiler simply maps the problem
specification onto three proximal operators
\begin{verbatim}
objective:
  add(
    least_squares(add(dense(C)*var(x), scalar(-1.00)*const(d))),
    norm1(var(y)),
    tv_1d(var(z)))

constraints:
  zero(add(var(z), scalar(-1.00)*var(y)))
  zero(add(var(z), scalar(-1.00)*var(x)))
\end{verbatim}
with the addition of the necessary variable copies and equality constraints to
make the objective separable.

\begin{figure}
  \centering
  \includegraphics[scale=0.5]{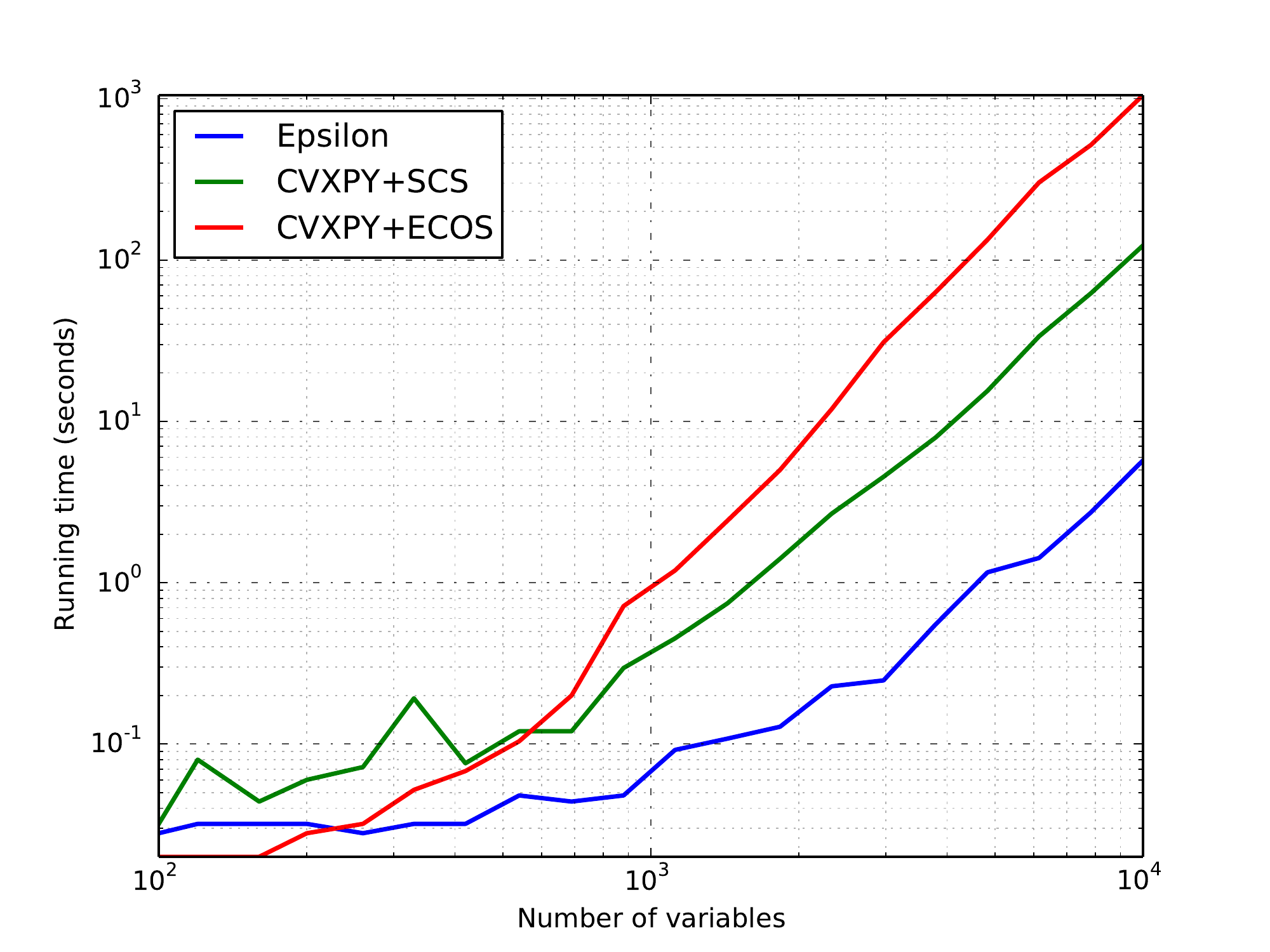}
  \caption{Comparison of running times on total variation example with $X \in
    \mathbb{R}^{m \times 10m}$, $y = X\theta_0 + \epsilon$ with $\epsilon
    \sim \mathcal{N}(0, 0.05^2)$ where $\theta_0$ is piecewise constant with
    segments of length 10.}
  \label{fig-scale-fused-lasso}
\end{figure}

Figure \ref{fig-scale-fused-lasso} compares Epsilon to the conic solvers on a
sequence of problems with $X \in \mathbb{R}^{m \times 10m}$. The main difference
between the two approaches is
that while Epsilon directly solves the proximal operator for the total variation
penalty (using the dynamic programming algorithm of \cite{johnson2013dynamic}),
while transforming to conic form requires reformulating the fused lasso penalty
as a linear program using auxiliary variables and involving the
finite differencing operator $D \in \{-1,0,1\}^{(n-1)
  \times n}$
\begin{equation}
D = \left [ \begin{array}{rrrrr}
1 & -1 & 0 & \cdots \\
0 & 1 & -1 & \cdots \\
0 & 0 & 1 & \cdots \\
\vdots & \vdots & \vdots & \ddots \\ \end{array} \right ]
\end{equation}
containing $1$ on the diagonal, $-1$ on the first super diagonal and 0 elsewhere.
For even moderate $n$ this linear operator (which corresponds to the edge
incidence matrix for the chain graph) is poorly conditioned which can be
problematic for general solvers, see e.g. \cite{ramdas2014fast} for further details. The
dedicated proximal operator avoids these issues, reducing running times---on a
problem with $10^4$ variables, Epsilon requires $5.7$ seconds vs. $123$ seconds for SCS.











\subsection{Library of convex programming examples}


In this final section we present results on a library of example convex problems
from statistical machine learning appearing frequently in the literature
(e.g. \cite{boyd2011distributed,o2013operator}). In general, each example
depends on a variety of factors such as dimensions of the constants,
data generation scheme and setting of the hyperparameters---we have chosen
reasonable defaults for these settings. The complete specification for each
problem is available in the the submodule \texttt{epsilon.problems}, distributed
with the Epsilon code base. On each example, we run each solver considered using the default
tolerances which for Epsilon and SCS correspond to moderate accuracy and high
accuracy for ECOS\footnote{Modifying the tolerances for an interior point method
  does not materially affect the comparison due to the nature of the
  bottlenecks.}. In practice, we observe that all solvers converge to a relative
accuracy of $10^{-2}$ which is reasonable for the statistical applications under
consideration.

\begin{table}
\centering
  \begin{tabular}{l|rr|rr|rr}
               &\multicolumn{2}{c|}{            Epsilon }&\multicolumn{2}{c|}{
      CVXPY+SCS }&\multicolumn{2}{c}{         CVXPY+ECOS } \\
 Problem       &   Time & Objective &   Time & Objective &   Time & Objective  \\
 \hline
 \texttt{basis\_pursuit}&   1.35s&$   1.44 \times 10^{2}$&  17.26s&$   1.45 \times 10^{2}$& 217.68s&$   1.45 \times 10^{2}$ \\
\texttt{covsel}&   0.93s&$   3.74 \times 10^{2}$&  25.09s&$   3.73 \times 10^{2}$&   -&        - \\
\texttt{fused\_lasso}&   3.87s&$   7.46 \times 10^{1}$&  57.85s&$   7.41 \times 10^{1}$& 641.34s&$   7.41 \times 10^{1}$ \\
\texttt{hinge\_l1}&   3.71s&$   1.15 \times 10^{3}$&  45.59s&$   1.15 \times 10^{3}$& 678.47s&$   1.15 \times 10^{3}$ \\
\texttt{hinge\_l1\_sparse}&  14.26s&$   1.38 \times 10^{3}$& 106.75s&$   1.38 \times 10^{3}$& 183.65s&$   1.38 \times 10^{3}$ \\
\texttt{hinge\_l2}&   3.58s&$   3.87 \times 10^{3}$& 133.10s&$   3.87 \times 10^{3}$&1708.31s&$   3.87 \times 10^{3}$ \\
\texttt{hinge\_l2\_sparse}&   1.82s&$   8.08 \times 10^{3}$&  28.40s&$   8.09 \times 10^{3}$&  47.72s&$   8.08 \times 10^{3}$ \\
\texttt{huber}&   0.20s&$   2.18 \times 10^{3}$&   3.17s&$   2.18 \times 10^{3}$&  28.43s&$   2.18 \times 10^{3}$ \\
\texttt{lasso}&   3.69s&$   3.21 \times 10^{1}$&  20.54s&$   3.21 \times 10^{1}$& 215.68s&$   3.21 \times 10^{1}$ \\
\texttt{lasso\_sparse}&  13.58s&$   4.37 \times 10^{2}$&  56.94s&$   4.37 \times 10^{2}$& 277.80s&$   4.37 \times 10^{2}$ \\
\texttt{least\_abs\_dev}&   0.10s&$   7.09 \times 10^{3}$&   2.96s&$   7.10 \times 10^{3}$&  11.46s&$   7.09 \times 10^{3}$ \\
\texttt{logreg\_l1}&   3.70s&$   8.18 \times 10^{2}$&  51.60s&$   8.18 \times 10^{2}$& 684.86s&$   8.17 \times 10^{2}$ \\
\texttt{logreg\_l1\_sparse}&   6.69s&$   9.61 \times 10^{2}$&  33.35s&$   9.63 \times 10^{2}$& 310.02s&$   9.61 \times 10^{2}$ \\
\texttt{lp}&   0.33s&$   7.77 \times 10^{2}$&   3.78s&$   7.75 \times 10^{2}$&   7.58s&$   7.77 \times 10^{2}$ \\
\texttt{mnist}&   0.91s&$   1.75 \times 10^{3}$& 219.63s&$   1.72 \times 10^{3}$&1752.97s&$   1.72 \times 10^{3}$ \\
\texttt{mv\_lasso}&   7.14s&$   4.87 \times 10^{2}$& 824.83s&$   4.88 \times 10^{2}$&    -&        - \\
\texttt{qp}&   1.39s&$   4.30 \times 10^{3}$&   3.20s&$   4.28 \times 10^{3}$&  23.12s&$   4.24 \times 10^{3}$ \\
\texttt{robust\_pca}&   0.59s&$   1.71 \times 10^{3}$&   2.88s&$   1.71 \times 10^{3}$&   -&        -\\
\texttt{tv\_1d}&   0.13s&$   2.29 \times 10^{5}$&  51.85s&$   2.95 \times 10^{5}$&  -&        - \\
\end{tabular}
\caption{Comparison of running time and objective value between Epsilon and
  CVXPY with SCS and ECOS solvers, a value of ``-'' indicates a lack of result
  either due to solver failure, unsupported problem or 1 hour timeout.}
\label{tab-comparison}
 \end{table}


The running times in Table \ref{tab-comparison} show that on all problem
examples considered, Epsilon is faster than SCS and ECOS and often by a wide
margin. In general, we observe that Epsilon tends to solve problems in fewer
ADMM iterations and for many problems the iterations are faster due in part to
operating on a smaller number of variables. There are numerous reasons for these
differences, some of which we have explored in the more detailed examples
appearing earlier in this section.

\section{Conclusions}

We have discussed Epsilon, a new system for general convex programming based on
representing any disciplined convex problem as a sum of functions with fast
proximal operators. The central idea is that by retaining more of the original
problem structure with a richer intermediate representation (the prox-affine
form), we achieve more efficient solution methods. This is inspired in part
by the recent surge of popularity in operator splitting methods leading to
a very large number of ``specialized'' algorithms which ultimately can be
reduced to a particular sequence of proximal operator invocations. In short,
Epsilon seeks to automate this process by providing a large library of proximal
(and linear) operators as well as a compiler to transform any DCP-valid
problem into a separable prox-affine form for direct solution by an ADMM-based
algorithm. Crucially, the iterative algorithm is agnostic to the details of the
operators and thus the system can be extended by adding new operator
implementations and the corresponding rules to the Epsilon compiler. Each new
proximal operator that is added to the Epsilon system in this fashion reduces
the need for laborious specialized implementations to solve a single problem. Ultimately,
an interesting direction for future work is thus the further development of Epsilon
and other tools toward automating many of the mechanical transformations and other
tasks required by convex programming practitioners to develop efficient
algorithms and deploy them on data at scale.

\bibliographystyle{plain}
\bibliography{epsilon}

\end{document}